\documentclass[12pt]{amsart}

\usepackage{amsmath,xspace,amssymb,mathrsfs}
\usepackage{color}

\input xy
\xyoption{all}
\xyoption{2cell}
\UseAllTwocells
\CompileMatrices

\newcommand{\Gr}{\operatorname{Gr}}

\renewcommand{\phi}{\varphi}

\newcommand{\dm}{\operatorname{dim}}

\newtheorem{proposition}{Proposition}[section]

\newtheorem{corollary}[proposition]{Corollary}
\newtheorem{theorem}[proposition]{Theorem}

\theoremstyle{definition}

\newtheorem{remark}[proposition]{Remark}

\begin{document}

\title{Grassmannians over a finite field}

\author[A. Tarizadeh]{Abolfazl Tarizadeh}
\address{Department of Mathematics, Faculty of Basic Sciences, University of Maragheh \\
P. O. Box 55136-553, Maragheh, Iran.
 }
\email{ebulfez1978@gmail.com}

\date{}

\begin{abstract} In this paper, we characterize the Grassmannian $\Gr(d,n)$ in terms of the row echelon forms of rank $d$. Using this characterization, then in the case of finite field we give a polynomial-type formula for the cardinality of the Grassmannian. \\
\end{abstract}

\maketitle

\section{Introduction}

If $A$ is an $m\times n$ matrix over a field $F$ then we denote by $\mathfrak{rs}(A)$ the subspace of $F^{n}$ generated by the row vectors of $A$ and it is called the row space of $A$. We denote the set of $n\times n$ row echelon forms over $F$ by $\mathscr{E}_{F}(n)$. Also, $\mathscr{E}_{F}(d,n)$ denotes the set of $n\times n$ row echelon forms over $F$ of rank $d$. For the definition of row echelon form (row-reduced echelon matrix) please consider \cite[Chap. 1, p. 11]{Hoffman}. \\

The set of $d$ dimensional subspaces of an $n$ dimensional vector space over a field $F$ is denoted by $\Gr_{F}(d,n)$ or shortly by $\Gr(d,n)$ if there is no confusion on the base field and it is called the Grassmannian. Grassmannians have rich mathematical structures with wide applications and have been studied in the literature over the years from algebraic, topological and geometric points of view, see e.g. \cite{Fedorov}, \cite{Ghorpade}, \cite{Kleiman} and \cite{Neretin}. \\

If $F$ is a finite field of size $q$ then the well known Gaussiann formula computes the cardinality of the Grassmannian $\Gr(d,n)$ as a rational formula:\\
$$|\Gr(d,n)|=
\frac{(q^{n}-1)(q^{n}-q)...(q^{n}-q^{d-1})}{(q^{d}-1)(q^{d}-q)...(q^{d}-q^{d-1})} .$$ \\

In this paper we give a new formula for the size of the Grassmannian. In fact, Theorem \ref{Theorem I} and Theorem \ref{formula} are the main contributions of this paper. Theorem \ref{Theorem I} establishes a one to one correspondence between the Grassmannian $\Gr_{F}(d,n)$ and $\mathscr{E}_{F}(d,n)$. Specially, if the base field is finite then we obtain a polynomial-type formula for the cardinality of the Grassmannian $\Gr(d,n)$, see Theorem \ref{formula}. \\

\section{Main Results}

\begin{theorem}\label{Theorem I} The map $R\rightsquigarrow \mathfrak{rs}(R)$ is a bijection from $\mathscr{E}_{F}(n)$ onto $\bigcup\limits_{d=0}^{n}\Gr(d,n)$. \\
\end{theorem}

{\bf Proof.} Without loss of generality, we may work with the vector space $F^{n}$. First we show that the above map is surjective. If $W\in\bigcup\limits_{d=0}^{n}\Gr(d,n)$ then $\dm W= d$ for some $0\leq d\leq n$. Let $\mathcal{B}=\{\alpha_{1},...,\alpha_{d}\}$ be an ordered basis of $W$ where $\alpha_{i}\in F^{n}$. Let $A$ be the $n\times n $ matrix over $F$ whose first $d$ rows are the vectors of $\mathcal{B}$ and the remaining $n-d$ rows of $A$ are zero. It is clear that $\mathfrak{rs}(A)=W$. It is well known that there exists an $n\times n$ row echelon form $R$ over $F$ which is row-equivalent to $A$. By \cite[Chap. 2, Theorem 9]{Hoffman}, row-equivalent matrices have the same row spaces. Thus $\mathfrak{rs}(R)=\mathfrak{rs}(A)=W$. Hence, the above map is surjective. Suppose $\mathfrak{rs}(R)=\mathfrak{rs}(R')$ for some $R, R'\in\mathscr{E}_{F}(n)$. We shall prove that $R=R'$. Let $\rho_{1},...,\rho_{d}$ (resp. $\rho'_{1},...,\rho'_{d'}$) be the non-zero row vectors of $R$ (resp. $R'$). Assume that the leading non-zero entry of $\rho_{i}$ (resp. $\rho'_{j}$) occurs in column $k_{i}$ (resp. $k'_{j}$). We have $\dm\big(\mathfrak{rs}(R)\big)=\dm\big(\mathfrak{rs}(R')\big)$. By \cite[Chap. 2, Theorem 10]{Hoffman},  $d=\dm\big(\mathfrak{rs}(R)\big)$ and $d'=\dm\big(\mathfrak{rs}(R')\big)$. Therefore $d=d'$. Now to prove $R=R'$ it suffices to show that $\rho_{i}=\rho'_{i}$ for all $i$ with $1\leq i\leq d$. If $\beta=(b_{1},...,b_{n})\in\mathfrak{rs}(R)$ then there exist scalers $c_{1},...,c_{d}\in F$ such that $\beta=\sum\limits_{i=1}^{d}c_{i}\rho_{i}$. We have $c_{\ell}=b_{k_{\ell}}$ for all $\ell$. Because if $\rho_{i}=(R_{i1},...,R_{in})$ then
from $\beta=(b_{1},...,b_{n})=\sum\limits_{i=1}^{d}c_{i}\rho_{i}$ we get that $b_{k_{\ell}}=\sum\limits_{i=1}^{d}c_{i}R_{i k_{\ell}}=\sum\limits_{i=1}^{d}c_{i}\delta_{i\ell}=c_{\ell}$. Therefore: \begin{equation}\label{equ I} \beta=\sum\limits_{i=1}^{d}b_{k_{i}}\rho_{i}. \end{equation} The expression (1), yields that if $\beta\neq 0$ then the index of the first non-zero component of $\beta$ belongs to the set $\{k_{1},...,k_{d}\}$. Because at least one of the $b_{k_{1}},...,b_{k_{d}}$
is non-zero. Suppose $b_{k_{t}}$ is the first non-zero between them, i.e. $k_{t}$ is the least index between the indices $\{k_{1},...,k_{d}\}$ for which $b_{k_{t}}\neq 0$. Therefore we can write: \begin{equation}\label{equ 2} \beta=\sum\limits_{i=t}^{d}b_{k_{i}}\rho_{i}. \end{equation}
We prove that $b_{j}=0$ for all $j$ with $j< k_{t}$. From (\ref{equ 2}), we have $b_{j}=\sum\limits_{i=t}^{d}b_{k_{i}}R_{ij}$. Since $j<k_{t}<k_{t+1}<...<k_{d}$, then by the definition of row-reduced echelon matrix, $R_{ij}=0$ for all $i$ with $t\leq i\leq d$. Hence, $b_{j}=0$. Therefore, $b_{k_{t}}$ is the first non-zero component of $\beta$. On the other hand, since $\rho'_{j}=(R'_{j1},...,R'_{jn})\in\mathfrak{rs}(R')=\mathfrak{rs}(R)$ thus the expression (\ref{equ I}) implies that $\rho'_{j}=\sum\limits_{i=1}^{d}R'_{jk_{i}}\rho_{i}$ for all $j$ with $1\leq j\leq d$. By the definition of row-reduced echelon matrix, the first non-zero component of $\rho'_{j}$ occurs in $k'_{j}$, i.e. $R'_{jk'_{j}}=1$. Now, using what we have just proved in the above for the first non-zero component, we have $k'_{j}\in\{k_{1},...,k_{d}\}$ for all $j$ with $1\leq j\leq d$. Also, since $k_{1}<k_{2}<...<k_{d}$ (resp. $k'_{1}<k'_{2}<...<k'_{d}$), one has then $k'_{j}=k_{j}$ for all $j$ with $1\leq j\leq d$. Finally, we have $\rho'_{j}=\sum\limits_{i=1}^{d}R'_{jk_{i}}\rho_{i}=
\sum\limits_{i=1}^{d}R'_{jk'_{i}}\rho_{i}=
\sum\limits_{i=1}^{d}\delta_{ji}\rho_{i}=\rho_{j}$ for all $j$. Therefore $R=R'$. $\Box$ \\

\begin{corollary}\label{Corollary I} There exists a one-to-one correspondence between the Grassmannian $\Gr_{F}(d,n)$ and $\mathscr{E}_{F}(d,n)$. \\
\end{corollary}

{\bf Proof.} It is an immediate consequence of Theorem \ref{Theorem I} and \cite[Chap. 2, Theorem 10]{Hoffman}. $\Box$ \\

In Theorem \ref{formula}, by $S^{(d)}$ we mean the set of all $(s_{1},...,s_{d})\in\mathbb{N}^{d}$ such that $1\leq s_{1}\leq n-d+1$ and
$s_{i-1}< s_{i}\leq n-d+i$ for all $i$ with $2\leq i\leq d$. \\

\begin{theorem}\label{formula} If $F$ is a finite field with $q$ elements Then:
  $$|\Gr_{F}(d,n)|=\sum\limits_{(s_{1},..., s_{d})\in S^{(d)}}q^{d(n-d)+\frac{d(d+1)}{2}-\sum\limits_{i=1}^{d}s_{i}}.$$ \\
\end{theorem}

{\bf Proof.} For each $\mathbf{s}=(s_{1},...,s_{d})\in S^{(d)}$ consider the subset $\mathcal{R}_{\mathbf{s}}\subseteq\mathscr{E}_{F}(d,n)$ consisting of all $R\in\mathscr{E}_{F}(d,n)$ such that for each $1\leq i\leq d$ the leading non-zero entry of $i-$th row of $R$ occurs in the column $s_{i}$. Clearly the $\mathcal{R}_{\mathbf{s}}$ are pairwise disjoint sets and
$\mathscr{E}_{F}(d,n)=\bigcup\limits_{\mathbf{s}\in S^{(d)}}\mathcal{R}_{\mathbf{s}}$. Therefore by Corollary \ref{Corollary I}, we have
$|\Gr(d,n)|=\sum\limits_{\mathbf{s}\in S^{(d)}}|\mathcal{R}_{\mathbf{s}}|$. It suffices to show that: $$|\mathcal{R}_{\mathbf{s}}|=
q^{d(n-d)+\frac{d(d+1)}{2}-\sum\limits_{i=1}^{d}s_{i}}.$$
If $R\in\mathcal{R}_{\mathbf{s}}$ then for each $1\leq i\leq d$,
at each $i-$th row of $R$, $d$ entries are $0$ or $1$ (the entry $1$ occurs in the column $s_{i}$ and zero entries in the columns $s_{j}$, $j\neq i$). Also, on every $i-$th row of $R$, by the definition of row echelon form, we have $R_{j, s_{i}}=0$ for each $j < s_{i}$, the number of these entries is $s_{i}-1$, since we have counted $i-1$ of them already in the above, therefore on every $i-$th row of $R$, the total number of entries which are $0$ or $1$ is $d+(s_{i}-1)-(i-1)=d+s_{i}-i$. Thus on every $i-$th row of $R$, $n-d-(s_{i}-i)$ entries are arbitrary scalers in $F$, and so in matrix $R$, $\sum\limits_{i=1}^{d}(n-d-(s_{i}-i))
=d(n-d)+\frac{d(d+1)}{2}-\sum\limits_{i=1}^{d}(s_{i})$ entries are arbitrary scalers. Therefore we get that $|\mathcal{R}_{\mathbf{s}}|=
q^{d(n-d)+\frac{d(d+1)}{2}-\sum\limits_{i=1}^{d}(s_{i})}$.  $\Box$ \\

\begin{remark} By Theorem \ref{formula}, we have  $|\Gr(1,n)|=\sum\limits_{i=0}^{n-1}q^{i}=q^{n-1}+q^{n-2}+...+1$. $|\Gr(0,n)|=|\Gr(n,n)|=1$. More generally, $|\Gr(d,n)|=|\Gr(n-d,n)|$. \\
\end{remark}

\begin{remark} In order to express the coefficients of the polynomial-type formula which appear in Theorem \ref{formula} more precisely we act as follows. We consider the equivalence relation $\sim$ on $S^{(d)}$ as
$(s_{1},...,s_{d})\sim(s'_{1},...,s'_{d})$ if and only if $\sum\limits_{j=1}^{d}s_{j}=\sum\limits_{j=1}^{d}s'_{j}$. For $(s_{1},...,s_{d})\in S^{(d)}$, let $[s_{1},...,s_{d}]$ denotes its equivalence class. Then Theorem \ref{formula} yields that: $$|\Gr(d,n)|=\sum\limits_{\ell=0}^{d(n-d)}c_{\ell}q^{d(n-d)-\ell}=
c_{0}q^{d(n-d)}+c_{1}q^{d(n-d)-1}+...+c_{d(n-d)}$$ where for each $0\leq\ell\leq d(n-d)$, $c_{\ell}$ is the cardinality of the class $[s_{1},...,s_{d}]$ whenever $\ell+\frac{d(d+1)}{2}=\sum\limits_{k=1}^{d}s_{k}$ for some $(s_{1},...,s_{d})\in S^{(d)}$.
By induction on $\ell$, one can show that each $c_{\ell}\geq 1$. For $\ell=0$, take $(s_{1}, s_{2},...,s_{d})=(1, 2,...,d)$ which belongs to $S^{(d)}$ and $0+\frac{d(d+1)}{2}=\sum\limits_{i=1}^{d}s_{i}$, thus $c_{0}=|[1,2,...,d]|\geq 1$. Now let $\ell>1$. Then by induction hypothesis, there exists $(s_{1},...,s_{d})\in S^{(d)}$ such that $\ell-1=\sum\limits_{i=1}^{d}(s_{i}-i)$. But there exists some $1\leq j\leq d$, so that $s_{j}<n-d+j$. Take $j_{0} :=\max\{j : 1\leq j\leq d, s_{j}< n-d+j\}$. Now, set $s'_{i}=s_{i}$ for each $i\neq j_{0}$ and $s'_{j_{0}}=s_{j_{0}}+1$. Then, $(s'_{1},...,s'_{d})\in S^{(d)}$ and $\ell=\sum\limits_{i=1}^{d}(s'_{i}-i)$. Therefore $c_{\ell}=|[s'_{1},...,s'_{d}]|\geq 1$. Furthermore, it is easy to see that $c_{0}=c_{1}=c_{d(n-d)-1}=c_{d(n-d)}=1$. Therefore $|\Gr(d,n)|=q^{d(n-d)}+q^{d(n-d)-1}+\mathbf{c_{2}}q^{d(n-d)-2}+
 \mathbf{c_{3}}q^{d(n-d)-3}...+\mathbf{c_{d(n-d)-2}}q^{2}+q+1$.  \\
We conclude this paper by proposing the following question. Could one describe the coefficients $\mathbf{c_{2}}, \mathbf{c_{3}},...,\mathbf{c_{d(n-d)-2}}$ which appear in the above formula somehow more precisely than described in the above? \\
\end{remark}


\begin{thebibliography}{10}
\bibitem{Fedorov}
R. Fedorov, Affine Grassmannians of group schemes and exotic principal bundles over affine line, Amer. J. Math. 138(4) (2016) 879-906.
\bibitem{Ghorpade}
S. R. Ghorpade and G. Lachaud, Hyperplane Sections of Grassmannians and the Number of MDS Linear Codes, Finite Fields and Their Applications, 7 (2001) 468-506.
\bibitem{Hoffman}
K. Hoffman and R. Kunze, Linear Algebra, 2nd ed., Prentice Hall, 1971.
\bibitem{Kleiman}
S. Kleiman, Geometry on Grassmannians and applications to splitting bundles and smoothing cycles, Pub. Math. de l'IHES, 36 (1969) 281-297.
\bibitem{Neretin}
Y. A. Neretin, The space $L^{2}$ on semi-infinite Grassmannian over finite field, Adv. Math. 250 (2014) 320-350.
\end{thebibliography}
\end{document}